\documentclass{article}
\usepackage{stmaryrd}
\usepackage{mathrsfs}
\usepackage[centertags]{amsmath}
\usepackage{amsfonts, dsfont}
\usepackage{amssymb}
\usepackage{amsthm}
\usepackage{color}

\input xy
\xyoption{all}

\theoremstyle{plain}
\newtheorem{thm}{Theorem}[subsection]
\newtheorem{cor}[thm]{Corollary}
\newtheorem{lem}[thm]{Lemma}
\newtheorem{prop}[thm]{Proposition}
\newtheorem*{thm2}{Theorem}

\theoremstyle{definition}
\newtheorem*{remark}{Remark}
\newtheorem*{ack}{Acknowledgments}

\newcommand{\bd}{\begin{defn}}
\newcommand{\ed}{\end{defn}}
\newcommand{\bl}{\begin{lem}}
\newcommand{\el}{\end{lem}}
\newcommand{\bp}{\begin{prop}}
\newcommand{\ep}{\end{prop}}
\newcommand{\bt}{\begin{thm}}
\newcommand{\et}{\end{thm}}
\newcommand{\bc}{\begin{cor}}
\newcommand{\ec}{\end{cor}}
\newcommand{\br}{\begin{remark}}
\newcommand{\er}{\end{remark}}
\newcommand{\bdi}{\begin{diagram}}
\newcommand{\edi}{\end{diagram}}
\newcommand{\beq}{\begin{equation}}
\newcommand{\eeq}{\end{equation}}
\newcommand{\ba}{\begin{array}}
\newcommand{\ea}{\end{array}}
\newcommand{\bpf}{\begin{proof}}
\newcommand{\epf}{\end{proof}}

\newcommand{\Z}{\mathds{Z}}
\newcommand{\Q}{\mathds{Q}}
\newcommand{\Zp}{\mathds{Z}_{p}}
\newcommand{\Qp}{\mathds{Q}_{p}}
\newcommand{\al}{\alpha}
\newcommand{\be}{\beta}
\newcommand{\Ga}{\Gamma}
\newcommand{\ga}{\gamma}

 \DeclareMathOperator{\Gal}{Gal}
\DeclareMathOperator{\Hom}{Hom} \DeclareMathOperator{\rank}{rank}

\newcommand{\X}{\mathcal{X}}
\newcommand{\Y}{\mathcal{Y}}
\newcommand{\p}{\mathfrak{p}}

\newcommand{\ot}{\otimes}

\newcommand{\im}{\mathrm{im}\,}

\newcommand{\lra}{\longrightarrow}

\newcommand{\ps}[1]{\llbracket #1 \rrbracket}

\begin{document}
\title{On the Iwasawa asymptotic class number formula for $\Zp^r\rtimes\Zp$-extensions}
 \author{Dingli Liang\footnote{School of Mathematics and Statistics,
 Wuhan University, Wuhan, Hubei, 430072, P.R.China.
  E-mail: \texttt{ldlqxx827@gmail.com}} \quad
  Meng Fai Lim\footnote{School of Mathematics and Statistics $\&$ Hubei Key Laboratory of Mathematical Sciences,
Central China Normal University, Wuhan, Hubei, 430079, P.R.China.
 E-mail: \texttt{limmf@mail.ccnu.edu.cn}} }
\date{}
\maketitle

\begin{abstract} \footnotesize
\noindent Let $p$ be an odd prime and $F_{\infty,\infty}$ a $p$-adic Lie extension of a number field $F$ with Galois group isomorphic to $\Zp^r\rtimes\Zp$, $r\geq 1$. Under certain assumptions, we prove an asymptotic formula for the growth of $p$-exponents of the class groups in the said $p$-adic Lie extension. This generalizes a previous result of Lei, where he establishes such a formula in the case $r=1$. An important and new ingredient towards extending Lei's result rests on an asymptotic formula for a finitely generated (not necessarily torsion) $\Zp\ps{\Zp^r}$-module which we will also establish in this paper. We then continue studying the growth of $p$-exponents of the class groups under more restrictive assumptions and show that there is an asymptotic formula in our noncommutative $p$-adic Lie extension analogous to a refined formula of Monsky (which is for the commutative extension) in a special case.

\medskip
\noindent Keywords and Phrases: Iwasawa asymptotic class number formula, $p$-adic Lie extension.

\smallskip
\noindent Mathematics Subject Classification 2010: 11R23, 11R29, 11R20.
\end{abstract}

\section{Introduction}

The starting point of this paper is the following Iwasawa's celebrated asymptotic class number formula (cf. \cite[Theorem 11]{Iw59}).

\begin{thm2}[Iwasawa]
Let $F_{\infty}$ be a $\Zp$-extension of a number field $F$. Denote by $F_n$ the intermediate subfield of $F_{\infty}$ with index $|F_n:F|=p^n$. Write $e_n$ for the $p$-exponent of the $p$-class group of $F_n$. Then there exist $\mu, \lambda$ and $\nu$ (independent of $n$) such that
\[ e_n = \mu p^n+\lambda n +\nu\]
for $n\gg 0$.
\end{thm2}

This result is the first of its kind which describes the precise growth of the $p$-exponent of the $p$-class group in an infinite tower of number fields. Subsequently, it is natural to search for such a formula, either exact or asymptotic, for other infinite towers. This was established by Cucuo and Monsky in the situation when $F_{\infty}$ is a $\Zp^d$-extension of $F$ (cf. \cite[Theorem I]{CM}; also see \cite[Theorem 3.13]{Mon} for a more refined estimate, and it is this latter form which we state here).

\begin{thm2}[Cucuo-Monsky, Monsky]
Let $F_{\infty}$ be a $\Zp^d$-extension of a number field $F$. Here $e_n$ will denote the $p$-exponent of the $p$-class group of $F_n$, where $F_n$ is now the intermediate subfield of $F_{\infty}$ with $\Gal(F_n/F)\cong(\Z/p^n)^{d}$. Then there exist integers $\mu$ and $\lambda$, and a real number $\al^*$ (independent of $n$) such that
\[ e_n = \mu p^{dn}+\lambda np^{(d-1)n} + \al^*p^{(d-1)n}+ O(np^{(d-2)n}).\]
\end{thm2}

In view of the above mentioned results, the next natural direction of investigation is to consider the case of a noncommutative $p$-adic Lie extension and search for an analogue of an asymptotic formula. A stumbling block towards deriving such a formula for the noncommutative situation is that we do not have a nice enough structure theory for modules over noncommutative Iwasawa algebras unlike the commutative situation (see \cite{CFKSV, CSSalg}). Despite this, the intensive study on the Iwasawa theory of noncommutative $p$-adic Lie extension undertaken by Venjakob (see \cite{V02, V03}) gave him the impetus to predict that an analogous asymptotic formula should hold for a false-Tate extension (see \cite[pp. 187]{V03}).

To the best knowledge of the authors, the first result in this noncommutative direction was by Perbet \cite{Per}, where he obtained an asymptotic formula of the growth of the $p$-exponent of the $p^n$-torsion subgroup of the class group of the intermediate subfield $F_n$. However, Perbet's methods do not allow one to obtain a formula for the full $p$-class group. It was only very recently that Lei was able to make progress on this, where he obtained an asymptotic class number formula for a $\Zp\rtimes \Zp$-extension under certain hypothesis (see \cite[Corollary 5.3]{Lei} and below). In particular, Lei's result gave an affirmative answer to Venjakob's prediction.

The main goal of this paper is to extend Lei's result. Namely, we obtain an asymptotic class number formula for a $\Zp^{r}\rtimes \Zp$-extension under an analogous hypothesis to that of Lei (see Theorem \ref{main theorem}) which we now describe. To do this, we need to introduce some notation. Let $p$ be an odd prime and $r$ a positive integer. Fix a number field $F$. Let $F_{\infty,\infty}$ be a $p$-adic Lie extension of $F$ with Galois group $G=\Gal(F_{\infty,\infty}/F)$. The dimension of $G$ is then denoted by $d$ which is always assumed to be $\geq 2$. Suppose that $G$ contains a closed normal subgroup $H\cong\Zp^{r}$ such that $\Ga:=G/H \cong \Zp$. (Therefore, we have $d= r+1$.) Write $F^c$ for the subextension of $F_{\infty,\infty}$ fixed by $H$. The following ramification conditions $(R1)-(R4)$ are always assumed to be satisfied for our extension $F_{\infty,\infty}/F$.

$(R1)$ The number field $F$ has only one prime above $p$ which we denote by $\p$.

$(R2)$ The prime $\p$ is totally ramified in $F_{\infty,\infty}/F$.

$(R3)$ The set of primes of $F$ ramified in $F_{\infty,\infty}$ is finite and we shall denote this set by $\Sigma$.

$(R4)$ Every prime of $\Sigma$ is finitely decomposed in $F^c/F$.

For $0\leq m,n\leq \infty$, we let $H_m = H^{p^m}$, $\Ga_n = \Ga^{p^n}$ and $G_{n,m} = H_m\rtimes \Ga_n$. Set $F_{n,m}$ to be the fixed field of $G_{n,m}$. The $p$-exponent of the $p$-class group of $F_{n,m}$ is then denoted by $e_{n,m}$. Finally, we write $\X=\Gal(\mathcal{M}/F_{\infty,\infty})$, where $\mathcal{M}$ is the maximal unramified abelian pro-$p$ extension of $F_{\infty,\infty}$. The following is the main theorem of this paper noting that our $r$ here is $d-1$ of the previous authors.

\begin{thm2}[Theorem \ref{main theorem}] Retain the setting of the preceding paragraphs.
 Suppose further that $\X$ is finitely generated over $\Zp\ps{H}$. Then we have
\[ e_{n,n} = \rank_{\Zp\ps{H}}(\X)np^{rn} +O(p^{rn}).\]
\end{thm2}

We now say a little on the ideas of the proof. In \cite{Lei}, Lei established this formula for the case $r=1$. There he made use of the analysis of torsion $\Zp\ps{H}$-modules in Cucuo and Monsky's work \cite{CM}. Note that the hypothesis of the theorem only assumes that $\X$ is finitely generated over $\Zp\ps{H}$. In the situation of $r=1$, there is a nice enough structure theory for $\Zp\ps{H}$-modules, and combined with the fact that $H$ is pro-cyclic (since $r=1$), Lei was able to reduce the problem to considering torsion $\Zp\ps{H}$-modules which in turn allowed him to obtain a crucial estimate (see \cite[Proposition 5.2]{Lei}) required for his eventual proof. Unfortunately, when $r\geq 2$, the structure theory for $\Zp\ps{H}$-modules is less refined, and more importantly, the group $H$ is no longer pro-cyclic. Therefore, Lei's approach does not seem to be able to carry over. Hence this suggests the necessitation to work directly with finitely generated $\Zp\ps{H}$-modules which are not necessarily torsion. This is precisely the core of Section \ref{algebra}, where we obtain the following asymptotic formula for a finitely generated (not necessarily torsion) $\Zp\ps{H}$-modules analogous to those in \cite{CM}.

\begin{thm2}[Theorem \ref{general estimate}]
 Let $H=\Zp^r$, $r\geq 1$. Let $M$ be a finitely generated (not necessarily torsion) $\Zp\ps{H}$-module. Then we have \[e(M_{H_m}) = \mu_H(M)p^{rm}+O(mp^{(r-1)m}).\]
\end{thm2}

Here $e(M_{H_m})$ denotes the $p$-exponent of the torsion subgroup of $M_{H_m}$ and $\mu_H(M)$ denotes the $\mu_H$-invariant of the module $M$ (see Section 2 for the definition). Our estimate is cruder than those in \cite{CM, Mon} which is somewhat expected, since we are now working with possibly nontorsion modules. Fortunately, this crude estimate is enough for us to establish the analogue of the crucial estimate of Lei for $r\geq 2$ (see Proposition \ref{Lei crucial estimate}) which in turn allows us to prove our main theorem.

The remainder of the paper is concerned with (raising) questions on what we can say about the asymptotic growth under more stringent hypothesis on the structure of the Galois group $\X$. The motivation behind this study stems from the refined asymptotic formula of Monsky \cite[Theorem 3.13]{Mon} (or see above), where he obtained an asymptotic formula up to an error term of $O(np^{(d-2)n})$. It is then natural to ask if our estimates can be improved in this direction for the class of noncommutative $p$-adic Lie extensions considered in this paper. Unfortunately, we do not see any way of doing this in general at this point of writing which is why we took the following more restricted approach of study. Namely, we would like to ask for a possible shape of the asymptotic formula when $\X$ is a finitely generated \textit{torsion} $\Zp\ps{H}$-module? We remark that this line of approach is also inspired by Lei's and our approach in obtaining the asymptotic formulas, where we investigate the growth under an extra assumption on $\X$ which essentially eliminates the leading term $p^{(r+1)n}$ (or $p^{dn}$ in Cucuo-Monsky's notation). Therefore, in this more restrictive context that we adopted, our question is as follows.

\smallskip
Suppose that $\X$ is a finitely generated torsion $\Zp\ps{H}$-module, does there exist some $\al^*$ such that
\[e_{n,n} = \al^*p^{rn}+ O(np^{(r-1)n})?\]

\noindent
Note that it follows from our Theorem \ref{main theorem} and the assumption in the question that we have $e_{n,n} = O(p^{rn})$ (see Corollary \ref{main theorem corollary}) and hence the above question makes sense. Had $G_{\infty,\infty}$ been abelian, the result of Monsky (\cite[Theorem 3.13]{Mon}; also see beginning of this section) will confirm this speculation (Again, we note that our $r$ is Monsky's $d-1$).
At present, the best we can do in our noncommutative setting is the following estimate.

\begin{thm2}[Theorem \ref{main theorem3}]
 Suppose that $\X$ is a finitely generated torsion $\Zp\ps{H}$-module. Then we have
\[ e_{n,n} \leq \be^* p^{rn}+ O(np^{(r-1)n})\]
for some nonnegative integer  $\be^*$.
\end{thm2}

In fact, the $\be^*$ appearing in our upper bound has a description of the sum of certain $\mu_H$-invariants (see the proof of Theorem \ref{main theorem3}). However, we have to confess that, at this point of writing, we do not know whether there is any relation between our $\be^*$ and Monsky's $\al^*$.  Despite so, our description of $\be^*$ enables us to obtain the following estimate by imposing an even stronger assumption on $\X$, thus answering our question in a special case.

\begin{thm2}[Theorem \ref{main theorem2}]
 Suppose that $\X$ is a finitely generated torsion $\Zp\ps{H}$-module with trivial $\mu_H$-invariant. Then we have
\[ e_{n,n} = O(np^{(r-1)n}).\]
\end{thm2}

We then also discuss some interesting examples
(see Subsection 3.3). In particular, in one of the examples, we specialize to the extension $F_{\infty, \infty} = \Q(\mu_{p^{\infty}}, p^{-p^{\infty}})$ and show that $e_{n,n}=O(n)$ for an irregular prime $p<1000$ by combining the above theorem with results of Sharifi \cite{Sh}. This estimate, together with another result of Cucuo-Monsky \cite[Theorem II]{CM} in the commutative situation, naturally leads us to ask whether there exist $\lambda$ and $\nu$ such that $e_{n,n} = \lambda n +\nu$ for sufficiently large $n$? Unfortunately, we (again!) do not have an answer to this at this point of writing.

\begin{ack}
The authors like to thank Xin Wan for arranging the collaboration of the authors. M. F. Lim would also like to thank Antonio Lei for many insightful discussion on his paper \cite{Lei} and the subject on the asymptotic class number formula. The authors are also grateful to Antonio Lei for answering many of their questions and explaining certain technical aspects of his paper during the preparation of this paper. The authors also like to thank Sohei Tateno for pointing out an inaccuracy in an earlier version of the paper.
Some part of the research of this article was conducted when M. F. Lim was visiting the National Center for Theoretical Sciences of Taiwan during Feb 2018, and he would like to acknowledge the hospitality
and conducive working conditions provided by the said institute. This work contains material which forms part of D. Liang's undergraduate thesis. Finally, M. F. Lim's research is supported by the
National Natural Science Foundation of China under Grant No. 11550110172 and Grant No. 11771164.
 \end{ack}

\section{Algebraic preliminaries} \label{algebra}

As before, $p$ will denote a fixed odd prime. Throughout this section, $H$ will always denote a multiplicative group which is isomorphic to a direct sum of $r$ copies of the additive group of the $p$-adic ring of integers, where $r\geq 1$. For each $m$, write $H_m = H^{p^m}$. If $N$ is a $\Zp$-module, denote by $N(p)$ the submodule of $N$ consisting of elements of $N$ which are annihilated by some power of $p$. In the event that $N$ is finitely generated over $\Zp$, we write $e(N)$ for the $p$-exponent of $N(p)$, i.e., $|N(p)| = p^{e(N)}$.

\subsection{Some technical lemmas}
In this subsection, we collect several technical results required for subsequent discussion. We first make a remark that if $M$ is a finitely generated $\Zp\ps{H}$-module, then $H_i(H_m,M)$ is finitely generated over $\Zp$ (for instance, see \cite[Lemma 3.2.3]{LS}). Thus, it makes sense to speak of $\rank_{\Zp}H_i(H_m,M)$ and $e(H_i(H_m,M))$ which we do without further comment. We begin by showing that the exponent of the $\Zp$-torsion subgroup of $H_i(H_m,M)$ is bounded by a linear function.

\bl \label{tech lemma}
 Let $M$ be a finitely generated $\Zp\ps{H}$-module. Then there exists $c$ (independent of $m$) such that $p^{rm+c}$ annihilates $H_i(H_m,M)(p)$ for every $m$ and $i$.
\el

\bpf
 Since $\Zp\ps{H}$ has finite global dimension $r+1$ (cf. \cite[Page 288, Exercise 5]{NSW}), we have a finite free resolution
 \[ 0\lra P_{r+1}\lra \cdots \lra P_{1}\lra P_{0}\lra M\lra 0\]
 of $M$. Set $K_i=\ker(P_{i}\lra P_{i-1})$, where $P_{-1}$ is understood to be $M$. Since there are only finitely many $K_i$'s, we may apply \cite[Theorem 2.8]{CM} to find a $c$ (independent of $m$) such that $p^{rm+c}$ annihilates $K_i(p)$ for every $i$. A straightforward exercise in homological algebra tells us that $H_i(H_m,M)\hookrightarrow (K_{i-1})_{H_m}$ which in turn implies that $p^{rm+c}$ annihilates $H_i(H_m,M)(p)$ for every $m$ and $i$.
\epf

Before stating the next result, we recall the definition of a structure in the sense of Cucuo-Monsky \cite{CM}.
Let $M$ be a finitely generated $\Zp\ps{H}$-module. A structure $\mathcal{S}$ on $M$ consists of an integer $m_0$ and a finite set of pairs $(\tau_i, M_i)$ where $\tau\in H-H_1$ and the $M_i$ are submodules of $M$ (cf. \cite[Definition 4.1]{CM}). For $m\geq m_0$ and for each $i$, set
\[ \al_{i,m,m_0} =  \frac{\tau_i^{p^m}-1}{\tau_i^{p^{m_0}}-1}.\]
We then define $A_m(\mathcal{S},M)$ to be the submodule
\[ I_{H_n}M +\sum_i \al_{i,m,m_0}M_i\]
of $M$, where here $I_{H_n}$ denotes the augmentation ideal of $H_n$ in $\Zp\ps{H}$. We can now state the following analogue of Lemma \ref{tech lemma} for structures.

\bl \label{structure tech lemma}
 Let $M$ be a finitely generated $\Zp\ps{H}$-module and $\mathcal{S}$ a structure on $M$. Then there exists $c$ (independent of $m$) such that $p^{(r+1)m+c}$ annihilates the torsion subgroup of $M/A_m(\mathcal{S},M)$ for every $m\geq m_0$.
\el

\bpf
This is proven in \cite[Theorem 4.5]{CM} for a torsion $\Zp\ps{H}$-module but it can be checked that the same proof carries over for general modules.
\epf

We record another technical lemma which has been utilized in \cite{CM} without mention. For clarity of presentation, we have stated this and its corollary which will be frequently inferred in subsequent discussion of this paper.

\bl \label{short exact Zp}
 Let $0\lra A\stackrel{i}{\lra} B\stackrel{g}{\lra} C\lra 0$ be a short exact sequence of finitely generated $\Zp$-modules. Suppose that $p^j$ annihilates $A(p), B(p)$ and $C(p)$. Then we have an exact sequence
 \[ 0\lra A(p)\lra B(p) \lra C(p) \lra A_f/p^j\] of finite abelian groups,
 where here $A_f = A/A(p)$. Consequently, we have the following inequalities
 \[ |e(A)-e(B)| \leq  e(C),\] \[ |e(B)-e(C)| \leq  e(A) + j\rank_{\Zp}(A),\]
 \[ e(C) \leq  e(B)  + j\rank_{\Zp}(A).\]
\el

\bpf
For a finitely generated $\Zp$-module $N$, we have an identification $N = N(p)\oplus N_f$, where $N_f = N/N(p)$. Let $x\in N$. We then write $x = (x_t, x_f)$ which corresponds to such an identification, where $x_t\in N(p)$ and $x_f\in N_f$.
An application of the snake lemma to the diagram
 \[ \SelectTips{eu}{}
\xymatrix{
  0 \ar[r] & A \ar[d]^{p^j} \ar[r]^{} & B \ar[d]^{p^j}
  \ar[r]^{} & C \ar[d]^{p^j} \ar[r]& 0 \\
  0 \ar[r] & A \ar[r] & B
  \ar[r] & C \ar[r] &0   }
\]
yields an exact sequence
 \[ 0\lra A(p)\lra B(p) \lra C(p) \lra A/p^j.\]
It therefore remains to show that $\im(C(p) \lra A/p^f)\subseteq A_f/p^j$ which essentially amounts to a finer analysis of the snake lemma argument. Let $c = (c_t, 0)\in C(p)$. Then there exists $b=(b_t, b_f)\in B$ such that $g(b) = c$. It follows from our hypothesis on $p^j$ that $p^jb = (0, p^jb_f)$. From the snake lemma argument, there is a unique $a=(a_t, a_f)\in A$ being sent to $p^jb$ under the map $i$. Since the torsion component of $p^jb$ is trivial and $i$ is injective, we must have $a_t=0$. This in turn tells us that $c$ is sent to the class of $(0,a_f)$ in $A/p^j$ and hence lies in $A_f/p^j$.
\epf

\bc \label{four terms}
Let $A\lra B\lra C\lra D$ be an exact sequence of finitely generated $\Zp$-modules. Suppose that $p^j$ annihilates $A(p), B(p), C(p)$ and $D(p)$. Then we have
\[ |e(B)-e(C)| \leq  e(A) +e(D) + j\rank_{\Zp}(A).\]
\ec

\bpf
 Write $U=\ker (A\lra B)$, $V= \ker(B\lra C)$ and $W=\ker(C\lra D)$ and $D'=\mathrm{im}(C\lra D)$. It can be readily verified that $U(p),V(p), W(p)$ and $D'(p)$ are also annihilated by $p^j$. Therefore, we have
 \[ \ba{rl}
  |e(B)-e(C)| \!
 &\leq~  |e(B)-e(W)| + |e(W)-e(C)|\\
 &\leq~ e(V) + j\rank_{\Zp}(V) + e(D') \\
 &\leq~ e(A) +j\rank_{\Zp}(U) + j\rank_{\Zp}(V) + e(D) \\
 & = e(A) +j\rank_{\Zp}(A) + e(D), \\
 \ea\]
 where the second and third inequalities follow from repeated applications of Lemma \ref{short exact Zp} on the appropriate short exact sequences, and the final equality is a consequence of the simple observation that $\rank_{\Zp}(U) + \rank_{\Zp}(V) = \rank_{\Zp}(A)$.
\epf

\subsection{Estimates for pseudo-null modules}
In this subsection, we estimate the $\Zp$-rank and $p$-exponent of $H_i(H_m, M)$ for a pseudo-null $\Zp\ps{H}$-module $M$, where we recall that a finitely generated $\Zp\ps{H}$-module is said to be pseudo-null if its localization at every prime ideal of height one is trivial. For $M_{H_m}$, such estimates are worked out in \cite[Theorem 3.2]{CM}. For the $p$-exponent of $H_i(H_m, M)$ when $M$ is a $p$-torsion module, this is done in \cite[Theorem 2.5.1]{LimCMu} for general $p$-adic Lie groups. The following proposition can thus be viewed as a generalization of these previous results when $H\cong \Zp^r$.

\bp \label{psuedo-null estimate}
 Let $M$ be a pseudo-null $\Zp\ps{H}$-module. Then we have
 \[ \rank_{\Zp}(M_{H_m}) = O(p^{(r-2)m}), \quad e(M_{H_m}) = O(p^{(r-1)m}) \]
 and
 \[\rank_{\Zp}(H_i(H_m, M)) = O(p^{(r-1)m}),\quad e(H_i(H_m,M)) = O(mp^{(r-1)m})\]
 for $i\geq 1$.
\ep

\bpf
 As mentioned above, the estimates for $M_{H_m}$ are contained in \cite[Theorem 3.2]{CM}. It therefore remains to establish the estimates for the higher homology groups. We first make the following remark. For a finitely generated $\Zp\ps{H}$-module $M$, it is an easy exercise to verify that $M(p)$ is a $\Zp\ps{H}$-submodule of $M$. Since
the ring $\Zp\ps{G}$ is Noetherian, the module $M(p)$ is therefore also finitely
generated over $\Zp\ps{H}$. As a consequence, one can find an integer
$t$ such that $p^t$ annihilates $M(p)$. It then follows that for every $i$, $H_i(H_m, M(p))$ is a finitely generated $\Zp$ module annihilated by $p^t$. This in turn implies that $H_i(H_m, M(p))$ is finite.

 Now consider the short exact sequence
 \[ 0\lra M(p) \lra M \lra M_f\lra 0,\]
 where we write $M_f := M/M(p)$. From this short exact sequence, we have an exact sequence
 \[ H_i(H_m, M(p))\lra  H_i(H_m,M) \lra  H_i(H_m, M_f) \lra H_{i-1}(H_m, M(p)). \]
 Since $H_i(H_m,M(p))$ and $H_{i-1}(H_m,M(p))$ are finite by the remark in the preceding paragraph, we have $\rank_{\Zp}(H_i(H_m,M)) = \rank_{\Zp}(H_i(H_m,M_f))$ and
 \[ e(H_i(H_m, M))\leq e(H_i(H_m,M(p))) + e(H_i(H_m,M_f)),\]
 where here we emphasis that the latter inequality follows from the finiteness of $H_i(H_n,M(p))$. (Warning: note that for a sequence $A\lra B\lra C$ exact at $B$ but not necessarily short exact, we \textit{do not} have $e(B)\leq e(A) + e(C)$ in general.)
 Now by \cite[Theorem 2.5.1]{LimCMu}, we have $e(H_i(H_m,M(p))) = O(p^{(r-1)m})$. Therefore, we are essentially reduced to proving the estimates under the assumption that $M(p)=0$ which we will do for the remainder of the proof. We proceed first by establishing the estimate for $\rank_{\Zp}(H_i(H_m, M))$. The $H_m$-homology of the short exact sequence $0\lra M\lra M \lra M/p\lra 0$ yields an exact sequence
  \[ H_{i+1}(H_m, M/p)\lra H_i(H_m, M) \lra H_i(H_m, M) \lra H_i(H_m, M/p).\]
  Note that
  \[ \rank_{\Zp}H_i(H_m, M) = e(H_i(H_m,M)/p) - e(H_i(H_m,M)[p]).\]
  From the above exact sequence, we see that $e(H_i(H_m,M)/p) \leq  e(H_i(H_m,M/p))$ and $e(H_i(H_m,M)[p]) \leq  e(H_{i+1}(H_m,M/p))$, and these latter two quantities are $O(p^{(r-1)m})$ by another application of \cite[Theorem 2.5.1]{LimCMu}. Hence this gives the required estimate for $\rank_{\Zp}(H_i(H_m, M))$.

We now turn to estimating $e(H_i(H_m,M))$. 
By Lemma \ref{tech lemma}, there exists $c$ (independent of $m$) such that $p^{rm+c}$ annihilates $H_1(H_m,M)(p)$ and $H_2(H_m,M)(p)$ for every $m$. From the short exact sequence
\[ 0\lra M\lra M\lra M/p^{rm+c}\lra 0,\]
we obtain the following exact sequence
\[ H_2(H_m, M)\lra H_2(H_m, M/p^{rm+c})\lra H_1(H_m,M)(p)\lra 0, \]
where the surjectivity follows from the fact that since $p^{rm+c}$ annihilates $H_1(H_m,M)(p)$, one has the equality $H_1(H_m,M)[p^{rm+c}]= H_1(H_m,M)(p)$.
Note that $p^{rm+c}$ clearly annihilates $H_1(H_m, M/p^{rm+c})$. We may then
apply Corollary \ref{four terms} (taking $D=0$) to obtain the equality
\[ e(H_1(H_m,M)) \leq e(H_2(H_m,M/p^{rm+c})) +(rm+c)\rank_{\Zp}(H_2(H_m, M)). \]
By a similar argument to that in \cite[Corollary 2.4]{Per}, one can show that there exists a constant $C$ independent of $m$ such that \[e(H_2(H_m,M/p^{rm+c}))\leq C(rm+c)p^{(r-1)m}= O(mp^{(r-1)m}).\]
 On the other hand, the rank estimate established in the previous paragraph yields
\[(rm+c)\rank_{\Zp}(H_2(H_m, M)) = O(mp^{(r-1)m}).\]
 Putting these estimates together, we obtain $e(H_1(H_m,M)) = O(mp^{(r-1)m})$ which gives the required bound for the first homology group. The estimates for higher homology groups can be proved similarly.
 \epf

We end the subsection with the following result which compares the variation of the $p$-exponents of the $H_n$-invariants of two pseudo-isomorphic $\Zp\ps{H}$-modules. Recall that two finitely generated $\Zp\ps{H}$-modules are said to be pseudo-isomorphic if there is a $\Zp\ps{H}$-homomorphism $\varphi:M\lra N$ whose kernel and cokernel are pseudo-null $\Zp\ps{H}$-modules. Note that the property of being pseudo-isomorphic is not a symmetric relation for nontorsion modules.

\bp \label{pseudo-null compare}
 Let $M$ and $N$ be two finitely generated (not necessarily torsion) $\Zp\ps{H}$-modules. Suppose that there is a $\Zp\ps{H}$-homomorphism $\varphi:M\lra N$ whose kernel and cokernel are pseudo-null $\Zp\ps{H}$-modules. Then for each $i\geq 0$, we have \[|e(H_i(H_m, M))-e(H_i(H_m,N))| = O(mp^{(r-1)m}).\]
\ep

\bpf
 The statement will follow if it holds in the two special cases of exact
sequences
\[\ba{c} 0\lra P \lra M \lra N\lra 0, \\
  0\lra M \lra N \lra P\lra 0, \ea
\] where $P$ is a finitely generated pseudo-null $\Zp\ps{H}$-module. We will prove the second case, the
first case has a similar argument. By Lemma \ref{tech lemma}, there exists $c$ which is independent of $m$ and such that $p^{rm+c}$ annihilates the torsion subgroups of $H_i(H_m, M), H_i(H_m, N)$ and $H_i(H_m,P)$ for every $i$ and $m$ (noting that $H_i(H_m,-)=0$ for $i\geq r+1$ and so there are only finite number of these groups, thus enabling one to find such a $c$). Applying $H_m$-invariant, we obtain an exact sequence
\[ H_{i+1}(H_m, P)\lra H_i(H_m,M)\lra H_i(H_m,N)\lra H_i(H_m,P). \]
It then follows from an application of Corollary \ref{four terms} that
\[|e(H_i(H_m,M)) - e(H_i(H_m,N))|\leq e(H_{i+1}(H_m,P))+ e(H_i(H_m,P)) + (rm+c)\rank_{\Zp}H_{i+1}(H_m,P).  \]
The required estimates now follow from Proposition \ref{psuedo-null estimate}.
\epf

\subsection{Estimates for elementary modules}

In this subsection, we prove the following estimate for modules of the form $M = \Zp\ps{H}/f^s$, where $f\in\Zp\ps{H}$ is a generator of a prime ideal of $\Zp\ps{H}$ of height one. For the next proposition, we write $\delta_{f,p} =1$ or $0$ accordingly as $f\Zp\ps{H}=p\Zp\ps{H}$ or $f\Zp\ps{H}\neq p\Zp\ps{H}$.

\bp \label{elementary estimate}
 Let $M = \Zp\ps{H}/f^s$, where $f\in\Zp\ps{H}$ is a generator of a prime ideal of $\Zp\ps{H}$ of height one. Then $\rank_{\Zp}(M_{H_m}) = \rank_{\Zp}(H_1(H_m,M))= O(p^{(r-1)m})$ and $e(M_{H_m}) = \delta_{f,p}sp^{rm}+O(mp^{(r-1)m})$.
 Furthermore, we have $H_i(H_m,M)=0$ for $i\geq 2$ and $e(H_i(H_m,M))=0$ for $i\geq 1$.
\ep

\bpf
 Since $M$ is torsion over $\Zp\ps{H}$, it follows from a formula of Harris (cf. \cite[Theorem 1.10]{Har00}) that $\rank_{\Zp}(M_{H_m}) = O(p^{(r-1)m})$. Now as
  $\Zp\ps{H}$ has no zero divisors, we have an exact sequence
 \[ 0\lra \Zp\ps{H} \stackrel{\cdot f^s}{\lra} \Zp\ps{H} \lra M \lra 0. \]
 Now note that since $\Zp\ps{H}$ is a free $\Zp\ps{H_m}$-module, we have $H_i(H_m, \Zp\ps{H}) = 0$ for $i\geq
1$. Therefore, by considering the $H_m$-homology, we obtain an
exact sequence
\[ 0\lra H_1(H_m, M)\lra \Zp[H/H_m] \lra \Zp[H/H_m] \lra M_{H_m} \lra 0 \]
and $H_i(H_m, M)=0$ for $i\geq 2$. The latter implies that $H_i(H_m,M)=0$ and $e(H_i(H_m,M))=0$ for $i\geq 2$. On other hand, it follows from the four terms exact sequence that $\rank_{\Zp}(M_{H_m}) = \rank_{\Zp}(H_1(H_m,M))$. Thus, this completes all the required estimates for $\Zp$-rank. Also, since $\Zp[H/H_m]$ has no $\Zp$-torsion, we have $e(H_1(H_m,M))=0$.

It therefore remains to establish the estimate for $e(M_{H_m})$. The proof proceeds as in \cite[Theorem 2.5]{CM} with some slight difference which we describe. Write $\zeta = (\zeta_1,...,\zeta_r)\in \mu_{p^{\infty}}^{\oplus r}$. Recall that two such $\zeta$ is said to be conjugate if there is an automorphism of $\bar{\Q}_p$ sending one to the other. By identifying $\Zp\ps{H}$ with $\Zp\ps{T_1,...,T_r}$, we may view every element $g$ of $\Zp\ps{H}$ as a power series in $r$ variables and it makes sense to speak of $g(\zeta-1): = g(\zeta_1-1,...\zeta_r-1)$ which is an element of $\Zp[\zeta]:=\Zp[\zeta_1,...,\zeta_r]$. Such an assignment gives rise to a homomorphism  $\Zp[H/H_m] \lra \oplus \Zp[\zeta]$. Summing these homomorphisms over the conjugacy class of $\zeta = (\zeta_1,...,\zeta_r)\in \mu_{p^n}$, we obtain a homomorphism $\varphi_n : \Zp[H/H_m] \lra \oplus \Zp[\zeta]$ which is called the cyclotomic embedding in \cite[Section 2]{CM}. By \cite[Theorem 2.2]{CM}, this is injective with a finite cokernel annihilated by $p^{rm}$. One can easily check that there is a commutative diagram
 \[\entrymodifiers={!! <0pt, .8ex>+} \SelectTips{eu}{}\xymatrix{
     \Zp[H/H_m] \ar[r] \ar[d]_{f} & \oplus \Zp[\zeta] \ar[d]_{f(\zeta -1)}\\
     \Zp[H/H_m] \ar[r] & \oplus \Zp[\zeta]  } \]
 where the vertical map are given by multiplication by the corresponding element labelled on the map. Note that the cokernel of the vertical map on the left is precisely $M_{H_m}$. Denote by $N_m$ the cokernel of the vertical map on the right. By \cite[Theorem 2.4 and Lemma 2.4]{CM}, we have
 \[ |e(M_{H_m}) - e(N_m)| \leq rm\rank_{\Zp}(\ker f)\]
 Since $\Zp[H/H_m]$ is finitely generated over $\Zp$, the latter is equal to $  rm\rank_{\Zp}(\mathrm{coker} f)= rm\rank_{\Zp}(M_{H_m})$ which is precisely $O(mp^{(r-1)m})$ by our rank estimate (this is the essential difference from that in \cite[Theorem 2.5]{CM}). It therefore remains to estimate $e(N_m)$. For this, one proceeds similarly to that in \cite[Theorem 2.5]{CM} and obtains $e(N_m) = \delta_{f,p}sp^{rm}+O(mp^{(r-1)m})$. The conclusion of the proposition is now immediate from these estimates.
\epf

\subsection{Estimates for general modules}

We now prove the estimates for a finitely generated (not necessarily torsion) $\Zp\ps{H}$-module. Before that, we recall the notion of the $\mu$-invariant which we do with slightly more generality. Let $\mathcal{G}$ be a uniform pro-$p$ group in the sense of \cite{DSMS}. It then follows from \cite[Proposition
1.11]{Ho2} (see also \cite[Theorem 3.40]{V02}) that there is a
$\Zp\ps{\mathcal{G}}$-homomorphism
\[ \varphi: M(p) \lra \bigoplus_{i=1}^t\Zp\ps{\mathcal{G}}/p^{\al_i},\] whose
kernel and cokernel are pseudo-null $\Zp\ps{\mathcal{G}}$-modules, and where
the integers $t$ and $\al_i$ are uniquely determined. The $\mu_{\mathcal{G}}$-invariant of $M$ is defined to be $\mu_{\mathcal{G}}(M) = \sum_{i=1}^t\al_i$. We can state the following main result of this section.

\bt \label{general estimate}
 Let $M$ be a finitely generated $\Zp\ps{H}$-module. Then $e(M_{H_m}) = \mu_H(M)p^{rm}+O(mp^{(r-1)m})$ and $e(H_i(H_m, M)) = O(mp^{(r-1)m})$ for every $i\geq 1$.
\et

The proof of the theorem will take up the remainder of this subsection. We begin by establishing the estimates for torsion $\Zp\ps{H}$-modules. We mention that even for torsion modules, our estimate is less precise than that in \cite[Theorem 3.4]{CM} due to the fact that we do not work under the stronger hypothesis on the rank growth as assumed there.

\bl \label{general estimate:torsion}
 Let $M$ be a finitely generated torsion $\Zp\ps{H}$-module. Then $e(M_{H_m}) = \mu_H(M)p^{rm}+O(mp^{(r-1)m})$ and
 $e(H_i(H_m, M)) = O(mp^{(r-1)m})$ for every $i\geq 1$.
\el

\bpf
By Proposition \ref{pseudo-null compare}, we may assume that $M$ is a direct sum of modules of the form $\Zp\ps{H}/f^s$, where $f$ is a generator of a prime ideal of height one. It then suffices to prove the result for each $\Zp\ps{H}/f^s$ which is precisely Proposition \ref{elementary estimate}.
\epf

We turn to the situation of a torsionfree $\Zp\ps{H}$-module. The following preparatory lemma is perhaps well-known to experts but for a lack of proper reference, we shall state this formally here.

\bl \label{embed fg}
 Let $M$ be a finitely generated $\Zp\ps{H}$-module of $\Zp\ps{H}$-rank $s$. Then there exists a map $M \lra \Zp\ps{H}^s$ with $\Zp\ps{H}$-torsion kernel and cokernel.

 In particular, if $M$ is a torsionfree $\Zp\ps{H}$-module, then $M$ injects into $\Zp\ps{H}^s$ with $\Zp\ps{H}$-torsion cokernel.
\el

\bpf
The proof is similar to that in the proof of \cite[Lemma 4.9]{LimFine}. For the convenience of the readers, we repeat the argument here. Let $K(H)$ denote the
field of fractions of $\Zp\ps{H}$. Write $M^+ = \Hom_{\Zp\ps{H}}(M, \Zp\ps{H})$. Then
by \cite[Proposition 2.5]{V02}, there is a canonical map $M\lra M^{++}$ with
$\Zp\ps{H}$-torsion kernel and cokernel. Choose $f_1,..., f_s\in
M^+$ such that they form a basis for $K(H)\ot_{\Zp\ps{H}}M^+$. Then
these elements give rise to a map $\Zp\ps{H}^s\lra M^+$ which clearly
has $\Zp\ps{H}$-torsion kernel and cokernel. Taking $\Zp\ps{H}$-dual, we
obtain a map $M^{++}\lra \Zp\ps{H}^s$ with $\Zp\ps{H}$-torsion kernel
and cokernel. Combining this with the above canonical map, we obtain
the required map. The second assertion is immediate from the first.
\epf

We can now prove our estimate for a torsionfree $\Zp\ps{H}$-module.

\bl \label{general estimate:torsionfree}
 Let $M$ be a finitely generated torsionfree $\Zp\ps{H}$-module with $\rank_{\Zp\ps{H}}(M)=s$. Then for each $i\geq 0$, we have $e(H_i(H_m,M))= O(mp^{(r-1)m})$.
\el

\bpf
By Lemma \ref{embed fg}, we have a short exact sequence
\[ 0\lra M\lra \Zp\ps{H}^s\lra N\lra 0,\]
for some torsion $\Zp\ps{H}$-module $N$. This short exact sequence in turn yields an exact sequence
\[ 0\lra H_1(H_m,N)\lra M_{H_m}\lra \Zp[H/H_m]^s\lra N_{H_m}\lra 0\]
and isomorphisms
\[ H_{i+1}(H_m,N)\cong H_i(H_m,M)\] for $i\geq 1$. Hence it follows that
$e(H_i(H_m,M))= e(H_{i+1}(H_m,N))$ for $i\geq 0$, where the equality for the case $i=0$ follows from the fact that $\Zp[H/H_m]^s$ has no $\Zp$-torsion. But since $N$ is torsion, the required estimate of the lemma is now a consequence of Lemma \ref{general estimate:torsion}.
\epf

We finally come to the proof of Theorem \ref{general estimate}.

\bpf[Proof of Theorem \ref{general estimate}]
By \cite[Proposition 5.1.7]{NSW}, there is a pseudo-isomorphism $M\stackrel{\sim}{\lra} T(M)\oplus M_{tf}$, where $T(M)$ is the $\Zp\ps{H}$-torsion submodule of $M$ and $M_{tf} = M/T(M)$. By Proposition \ref{pseudo-null compare}, we are thus reduced to proving the estimates for $T(M)$ and $M_{tf}$ which are precisely the contents of Lemmas \ref{general estimate:torsion} and \ref{general estimate:torsionfree}. Hence we have proven our theorem.
\epf

\subsection{An estimate for certain $\Zp\ps{\Ga}$-modules}

In this subsection, we quote the following useful result which will be used in the proof of our main theorem.

\bl \label{Gamma estimate}
Let $M$ be a finitely generated $\Zp\ps{\Ga}$-module with the properties that $M$ is finitely generated over $\Zp$ and that $M_{\Ga_n}$ is finite for every $n\geq 1$.
Let $f$ be a generator of the characteristic ideal of $M$. Let $n_0$ be an integer such that every irreducible distinguished polynomial that divides $f$ has degree $<p^{n_0-1}(p-1)$. Then
\[ e(M_{\Ga_n}) - e(M_{\Ga_0}) = \rank_{\Zp}(M)(n-n_0) + e(M(p)_{\Ga_n}) -e(M(p)_{\Ga_{n_0}}) .\]
for all $n\geq n_0$.   \el

\bpf
This is a special case of \cite[Proposition 4.6]{Lei} noting that $\mu_{\Ga}(M)=0$ as $M$ is finitely generated over $\Zp$.
\epf

\section{Arithmetic}

\subsection{Setup}
We recall the arithmetic setup following that in \cite{Lei}. As before, $p$ always denotes an odd prime and $r$ an integer $\geq 1$. Let $F$ be a number field. Denote by $F_{\infty,\infty}$ a $p$-adic Lie extension of $F$ whose Galois group $G=\Gal(F_{\infty,\infty})$ is isomorphic to $\Zp^r\rtimes \Zp$. Let $H$ be a closed normal subgroup of $G$ which is isomorphic to $\Zp^r$ and has the property that $\Ga:=G/H\cong\Zp$. We write $F^c$ for the fixed field of $H$. The ramification assumptions $(R1)-(R4)$ in the introductory section are always assumed to be valid.

For $0\leq m,n\leq \infty$, write $H_m = H^{p^m}$, $\Ga_n = \Ga^{p^n}$ and $G_{n,m} = H_m\rtimes \Ga_n$. The fixed field of $G_{n,m}$ is then denoted to be $F_{n,m}$. Let $\mathcal{M}$ be the maximal unramified abelian pro-$p$ extension of $F_{\infty,\infty}$. We write $\X=\Gal(\mathcal{M}/F_{\infty,\infty})$, $\Y=\Gal(\mathcal{M}/F)$ and $\Y_{n,m}=\Gal(\mathcal{M}/F_{n,m})$. Denote by $C_{n,m}$ the subgroup of $\Y_{n,m}$ generated by $[\Y_{n,m},\Y_{n,m}]$ and all the inertia groups in $\Y_{n,m}$ and write $B_{n,m} = C_{n,m}\cap\X$. The relation between $B_{n,m}$ and the class group of $F_{n,m}$ is established in \cite[Lemma 2.2 and Formula (2)]{Lei} which we record here.

\bp \label{exact sequence} We have $\X/B_{n,m}\cong Cl(F_{n,m})(p)$ and a short exact sequence
\[ 0\lra B_{n,m}/I_{G_{n,m}}\X\lra \X_{G_{n,m}}\lra Cl(F_{n,m})(p)\lra 0, \]
where $I_{G_{n,m}}$ is the augmentation ideal of $G_{n,m}$ in $\Zp\ps{G}$.
\ep

We also record the useful estimate of Lei \cite[Corollary 2.10]{Lei}.

\bp[Lei] \label{Lei}
$B_{n,m}/I_{G_{n,m}}\X$ is finitely generated over $\Zp$ and there exists $C$ independent of $1\leq n,m\leq \infty$ such that
\[\rank_{\Zp}(\X/B_{n,m})\leq Cp^{(r-1)m}\]
for every $m$ and $n$.
\ep

 From now on, we always write $e_{n,m}$ for $p$-exponent of $Cl(F_{n,m})(p)$. For each $n$, we also denote by $\X(\Ga_n)$ the $\Zp\ps{G}$-submodule of $\X$ which is generated by
 elements of the form $(\ga-1)x$, where $\ga\in \Ga_n$ and $x\in \X$.

\bp[Lei] \label{weak bound}
 For a fixed $n$, we have
 \[e_{n,m} = \mu_H(\X/\X(\Ga_n))p^{rm} + \rank_{\Zp\ps{H}}(\X/\X(\Ga_n))mp^{(r-1)m}+O(p^{(r-1)m}).\]
\ep

\bpf
 This is \cite[Corollary 3.4]{Lei} noting that our $r$ here is $d-1$ there.
\epf

\subsection{Main results}

We can now state and prove the main theorem of this paper which is a generalization of \cite[Corollary 5.3]{Lei}.

\bt \label{main theorem}
Retain the setting in the previous subsection. Suppose that $\X$ is
finitely generated over $\Zp\ps{H}$. Then we have
\[ e_{n,n} = \rank_{\Zp\ps{H}}(\X)np^{rn} +O(p^{rn}).\]
\et

Before proving our theorem, we first show the following crucial estimate which is established by Lei in the case $r=1$ (see \cite[Proposition 5.2]{Lei}). This estimate is now possible in general thanks to our analysis in Section \ref{algebra}.

\bp \label{Lei crucial estimate}
Retain all the assumptions of Theorem \ref{main theorem}. Then
\[ e(\X/B_{\infty,m}) \leq \mu_H(\X)p^{rm}+ O(mp^{(r-1)m}).\]
\ep

\bpf
Since $B_{\infty,m} =A_m(\mathcal{S},\X)$ for some structure $\mathcal{S}$ on $\X$ (cf. \cite[Proofs of Lemma 3.3 and Proposition 5.2]{Lei}), we may apply Lemmas \ref{tech lemma} and \ref{structure tech lemma} to obtain an $c$ which is independent of $m$ and has the property that $p^{(r+1)m+c}$ annihilates the torsion subgroups of $\X_{H_m}$ and $\X/B_{\infty,m}$ for every $m$. This in turn allows us to apply Lemma \ref{short exact Zp} to the exact sequence in Proposition \ref{exact sequence} which yields the following inequality
\[ e(\X/B_{\infty,m})\leq e(\X_{H_m}) + ((r+1)m+c)\rank_{\Zp}(B_{\infty,m}/I_{H_m}\X). \]
Thanks to Theorem \ref{general estimate}, we have that $e(\X_{H_m})=\mu_H(\X)p^{rm}+ O(mp^{(r-1)m})$. On the other hand, by virtue of Proposition \ref{Lei}, one has $((r+1)m+c)\rank_{\Zp}(B_{\infty,m}/I_{H_m}\X)=O(mp^{(r-1)m})$. The required bound of the proposition is now a consequence of these estimates.
\epf

We are in position to prove our theorem.

\bpf[Proof of Theorem \ref{main theorem}]
By \cite[Proposition 2.3]{CFKSV}, there exists a finite collection of $g_i\in\Zp\ps{G}$ such that each $\Zp\ps{G}/\Zp\ps{G}g_i$ is finitely generated over $\Zp\ps{H}$ and such that there is a surjection
\[ \bigoplus_i \Zp\ps{G}/\Zp\ps{G}g_i \twoheadrightarrow \X.\]
By \cite[Corollary 5.4]{DL}, there exists an integer $t$ (independent of $m$) such that the $\Zp\ps{\Ga}$-characteristic ideal of $\oplus_i \Zp\ps{G}/\Zp\ps{G}g_i$, and hence $\X_{H_m}$, factorizes into polynomials of degree $\leq t$. Since $\X/B_{\infty,m}$ is a quotient of $\X_{H_m}$, the same conclusion holds for $\X/B_{\infty,m}$.
Now fix a choice of $n_0$ such that $t< p^{n_0-1}(p-1)$. For now, fix an arbitrary positive integer $m$. As $(\X/B_{\infty,m})_{\Ga_n}\cong Cl(F_{n,m})(p)$ (cf. \cite[Lemma 2.2]{Lei}) is finite, we may apply Lemma \ref{Gamma estimate} to conclude that whenever $n\geq n_0$, one has
\[ |e_{n,m} - e_{n_0,m} - (n-n_0)\rank_{\Zp}(\X/B_{\infty,m})|= |e\big((\X/B_{\infty,m})(p)_{\Ga_n}\big)-e\big((\X/B_{\infty,m})(p)_{\Ga_{n_0}}\big) |\]
\[ \quad \quad \quad\leq e\big((\X/B_{\infty,m})(p)_{\Ga_n}\big)+e\big((\X/B_{\infty,m})(p)_{\Ga_{n_0}}\big) \leq 2e(\X/B_{\infty,m}). \]
Since $m$ is arbitrary, the above formula, in particular, holds when $m=n\geq n_0$. In other words, we have
\[ |e_{n,n} - e_{n_0,n} - (n-n_0)\rank_{\Zp}(\X/B_{\infty,n})|\leq 2e(\X/B_{\infty,n}) \]
for $n\geq n_0$. Therefore, in order to estimate $e_{n,n}$, it remains to estimate every other terms in this inequality.
As a start, Proposition \ref{Lei crucial estimate} tells us that $e(\X/B_{\infty,n})= O(p^{rn})$.
By Proposition \ref{weak bound}, we have \[e_{n_0,n} = \mu_H(\X/\X(\Ga_{n_{0}}))p^{rn} + O(np^{(r-1)n}) =   O(p^{rn})\] since $n_0$ is fixed.
On the other hand, it follows from \cite[Theorem 1.10]{Har00} that $\rank_{\Zp}(\X_{H_n})=
\rank_{\Zp\ps{H}}(\X)p^{rn}+O(p^{(r-1)n})$. Combining this with Propositions \ref{exact sequence} and \ref{Lei}, we obtain
\[\rank_{\Zp}(\X/B_{\infty,n}) = \rank_{\Zp\ps{H}}(\X)p^{rn}+O(p^{(r-1)n})\]
 which in turn implies that
\[(n-n_0)\rank_{\Zp}(\X/B_{\infty,n})= \rank_{\Zp\ps{H}}(\X)np^{rn}- \rank_{\Zp\ps{H}}(\X)n_0p^{rn}+O(np^{(r-1)n}) \]
\[ = \rank_{\Zp\ps{H}}(\X)np^{rn}+O(p^{rn}).\]
 The conclusion of the theorem now follows by combining these estimates.
\epf

As an immediate corollary, we have the following.

\bc \label{main theorem corollary}
Suppose that $\X$ is a
finitely generated torsion $\Zp\ps{H}$-module. Then we have
\[ e_{n,n} = O(p^{rn}).\]
\ec

As mentioned in the introductional section, we like to ask if there exists some $\al^*$ such that
\[e_{n,n} = \al^*p^{rn}+ O(np^{(r-1)n})\]
when $\X$ is a finitely generated torsion $\Zp\ps{H}$-module.
In the event that $G_{\infty,\infty}$ is commutative, this will follow from the result of Monsky (cf. \cite[Theorem 3.13]{Mon}). In our noncommutative situation, the best we can do at present is the following.

\bt \label{main theorem3}
Suppose that $\X$ is a
finitely generated torsion $\Zp\ps{H}$-module. Then we have
\[ e_{n,n} \leq \be^*p^{rn} + O(np^{(r-1)n})\]
for some $\be^*$.
\et

\bpf
By a similar argument to that in Theorem \ref{main theorem}, there exists a $n_0$ such that
\[ e_{n,n}\leq  e_{n_0,n} + \rank_{\Zp}(\X/B_{\infty,n})(n-n_0) +  2e(\X/B_{\infty,n}). \]
Performing the same calculations as before, we obtain the required estimate noting that this time round our calculation gives
\[ \be^* = \mu_H(\X) + \mu_{H}(\X/\X(\Ga_{n_0})). \]
\epf

As mentioned in the introduction, we are not able to relate our $\be^*$ with Monsky's $\al^*$. Despite so,
the calculations in the preceding proof gives us the following result which answers our question in a special case.

\bt \label{main theorem2}
Suppose that $\X$ is a
finitely generated torsion $\Zp\ps{H}$-module with $\mu_H(\X)=0$. Then we have
\[ e_{n,n} = O(np^{(r-1)n}).\]
\et

\bpf
As seen from the calculation of Theorem \ref{main theorem3}, we have
\[ e_{n,n} \leq (\mu_H(\X) + \mu_{H}(\X/\X(\Ga_{n_0}))p^{rn} + O(np^{(r-1)n})\]
for some $n_0$. By our hypothesis, we have $\mu_H(\X)=0$ which in turn implies that $\mu_{H}(\X/\X(\Ga_{n_0}))=0$. Hence the conclusion of corollary now follows.
\epf

\subsection{Examples}

We conclude the paper with some examples and further questions.

\medskip
(1) A nonzero integer $\al$ is said to be amenable for the prime $p$ if either $p|\al$ or $p \parallel \al^{p-1}-1$ holds. Denote by $\mu_p$ (resp., $\mu_{p^{\infty}}$) the multiplicative group of $p$-root of unity (resp., $p$-power root of unity). Let $F=\Q(\mu_p)$, $F^c= \Q(\mu_{p^{\infty}})$ and $F_{\infty, \infty} = \Q(\mu_{p^{\infty}}, \al_1^{-p^{\infty}},...,\al_r^{-p^{\infty}})$, where $\al_i$ are nonzero integers whose image in $\Qp^{\times}/(\Qp^{\times})^p$ are linearly independent over $\mathbb{F}_p$ and the products $\al_1^{n_1}\cdots\al_r^{n_r}$ are all amenable for $p$ for every $0\leq n_i\leq p-1$. It is well-known that the ramification properties $(R1), (R2)$ and $(R4)$ are satisfied. Property $(R3)$ is satisfied by \cite[Theorem 5.2 and Lemma 6.1]{Viv}. By the theorem of Ferrero-Washington \cite{FW} and a descent argument (for instance, see \cite[Proposition 2.1]{Sh}), we have that $\X$ is finitely generated over $\Zp\ps{H}$. Hence Theorem \ref{main theorem} applies in this situation.

In view of \cite[Question 1.3 and the two paragraphs after]{HS}, one would expect $\X$ to be a finitely generated torsion $\Zp\ps{H}$-module. Granted this, Corollary \ref{main theorem corollary} applies. However, it remains currently an open problem to prove this torsionness condition in general (but see below discussion).

\medskip
(2) Let $F=\Q(\mu_p)$, $F^c= \Q(\mu_{p^{\infty}})$ and $F_{\infty, \infty} = \Q(\mu_{p^{\infty}}, p^{-p^{\infty}})$. Now if $p$ is a regular prime, then the $p$-class class group of any intermediate field in this tower is trivial by a classical result of Iwasawa \cite{Iw56}. We therefore assume that the prime $p$ is irregular in our discussion here. This was the extension, where Venjakob \cite{V03} first made his guess on the form of the asymptotic class number formula and which we now know is correct thanks to the work of Lei \cite{Lei}.

Now if $p$ is an irregular prime $<1000$, a result of Sharifi \cite[Propositions 3.3 and 2.1a]{Sh} asserts that $\X$ is finitely generated over $\Zp$. In view of the results of Cucuo-Monsky for a commutative $p$-adic Lie extension (see \cite[Theorem II]{CM}), one is led to the following speculation:

\smallskip
Does there exist $\lambda$ and $\nu$ such that $e_{n,n} = \lambda n +\nu$ for sufficiently large $n$?
\smallskip

\noindent At this point of writing, we do not have an answer to this. We do however remark that it follows from Sharifi's result and our Theorem \ref{main theorem2} that one at least has the suggestive estimate $e_{n,n}=O(n)$.

For irregular primes $>1000$ and beyond, our current state of knowledge seems even less and we have nothing to say on this.

\end{document}